\documentclass[letterpaper, 10 pt, conference]{ieeeconf}

\IEEEoverridecommandlockouts
 \usepackage{amsmath, amssymb, bbm, xspace}
 \usepackage{amsfonts,mathrsfs}
 \usepackage{mathtools}
 \usepackage{color}
 \usepackage{graphicx}
 \usepackage{epsfig}
 \usepackage{algorithm}
 \usepackage{algorithmicx}
 \usepackage[acronym]{glossaries}
 \usepackage{subfigure}
 \usepackage{cancel}
 \usepackage{empheq}
 \usepackage{placeins}
 \let\labelindent\relax
 \usepackage{enumitem}
 \usepackage{textgreek}
 \usepackage{mathbbol}
 \usepackage{stackengine}
\usepackage{placeins}
\usepackage{dblfloatfix}

 \def\QEDhereeqn{\eqno\let\eqno\relax\let\leqno\relax\let\veqno\relax\hbox{\QED}}
 \def\QEDopenhereeqn{\eqno\let\eqno\relax\let\leqno\relax\let\veqno\relax\hbox{\QEDopen}}
 \newcommand{\bs}{\boldsymbol}
 \newcommand{\mc}{\mathcal}
 \renewcommand{\emph}{\textit}
 
 \newcommand{\0}{\bs 0}
 \def\1{{\bs 1}}
 \def\argmin{\mathop{\rm argmin}}
 \newcommand{\col}{\mathrm{col}}

 \newcommand{\proj}{\mathrm{proj}}
 \def\dom{\operatorname{dom}}
 \def\R{\mathbb{R}}
 \def\N{\mathbb{N}}
 \def\I{\mc{I}}
 \def\i{{i\in\mc{I}}}

 \DeclareSymbolFontAlphabet{\mathbbm}{bbold}
 \DeclareSymbolFontAlphabet{\mathbb}{AMSb}%
 
 \stackMath
	 \newcommand\tsup[2][2]{%
	 	\def\useanchorwidth{T}%
	 	\ifnum#1>1%
	 	\stackon[-.5pt]{\tsup[\numexpr#1-1\relax]{#2}}{\scriptscriptstyle\sim}%
	 	\else%
	 	\stackon[.5pt]{#2}{\scriptscriptstyle\sim}%
	 	\fi%
	 }

\makeglossaries
\newacronym{KKT}{KKT}{Karush--Kuhn--Tucker}
\newacronym{ADMM}{ADMM}{alternating direction method of multipliers}
\newacronym{OPF}{OPF}{optimal power flow}
\newacronym{OPFP}{OPFP}{optimal power flow problem}
\newacronym{NUM}{NUM}{network utility maximization}
\newacronym{LMI}{LMI}{linear matrix inequality}
\newacronym{BMI}{BMI}{bilinear matrix inequality}
\newacronym{LM}{LM}{Lyapunov-Metzler}
\newacronym{SDP}{SDP}{semidefinite programming}
\newacronym{LTI}{LTI}{linear time invariant}
\newacronym{MJLS}{MJLS}{Markov jump linear system}
\newacronym{PID}{PID}{proportional-integral-derivative}
\newacronym{PPA}{PPA}{proximal-point algorithm}
\newacronym{PPPA}{PPPA}{preconditioned proximal-point algorithm}
\newacronym{PPP}{PPP}{preconditioned proximal-point}
\newacronym{NE}{NE}{Nash equilibrium}

 \makeatletter \let\cl@part\relax \makeatother
 \usepackage{nohyperref}
 \usepackage[capitalize]{cleveref}
  \let\k\relax
 \def\k{{k \in \N}}  
 \usepackage{letltxmacro}

 \crefname{thm}{Theorem}{Theorems}
 \crefname{lem}{Lemma}{Lemmas}
 \crefname{cor}{Corollary}{Corollaries}
 \crefname{rem}{Remark}{Remarks}
 \crefname{alg}{Algorithm}{Algorithm}
 \crefname{figure}{Figure}{Figures}
 \crefname{assumption}{Assumption}{Assumptions}
 \crefname{corollary}{Corollary}{Corollaries}
\crefname{proposition}{Proposition}{Propositions}
\crefname{problem}{Problem}{Problems}
\crefname{example}{Example}{Examples}
\crefname{definition}{Definition}{Definition}
\crefname{condition}{C}{C}

 \crefname{thmlisti}{Theorem}{Theorem}
 \crefname{lemlisti}{Lemma}{Lemma}
 \crefname{asmlisti}{Assumption}{Assumptions}

 \newlist{thmlist}{enumerate}{1}
 \setlist[thmlist]{label=(\roman{thmlisti}), ref=\thethm(\roman{thmlisti}),noitemsep}
 \newlist{lemlist}{enumerate}{1}
 \setlist[lemlist]{label=(\roman{lemlisti}), ref=\thelem(\roman{lemlisti}),noitemsep}
 \newlist{asmlist}{enumerate}{1}
 \setlist[asmlist]{label=(\roman{asmlisti}), ref=\theassumption(\roman{asmlisti}),noitemsep,nosep,leftmargin=*} 

 \newtheorem{lemma}{Lemma}
 \newtheorem{theorem}{Theorem}
 \newtheorem{remark}{Remark}
 \newtheorem{assumption}{Assumption}

\newtheorem{proposition}{Proposition}
\newtheorem{example}{Example}
\newtheorem{definition}{Definition}

\newtheorem{condition}{C \hspace{-0.7em}} 

\usepackage{etoolbox}
\patchcmd{\smallmatrix}{\thickspace}{\kern.5em}{}{}

\usepackage{caption}
\usepackage{booktabs}
\usepackage{multirow}

\def\Fmc{\mc{F}}
\def\H{\mc{H}}
\def\zer{\operatorname{zer}}
\def\gra{\operatorname{gra}}
\def\Fbs{\bs{F}}
\def\Fa{\Fmc_{\alpha}} 
\def\Rmc{\mc{R}}
\def\Fop{\Rmc^\top \Fbs}
\def\condref#1{C\ref{#1}}
\def\subd#1{\partial_{#1}}
\def\grad#1{\nabla_{#1}}
\def\x{\bs{x}}
\def\sol{\mc{S}}
\def\X{\mc{X}}

\def\dom{\operatorname{dom}}
\def\A{\mc{A}_{\alpha}}
\def\Ak{\mc{A}_{\alpha^k}}
\def\Id{\operatorname{Id}}
\def\Wbs{\bs{W}}
\def\nc{\mathrm{N}}
\def\res{\mathrm{J}}
\def\ubs{\bs{u}}
\def\fix{\operatorname{fix}}
\def\y{\bs{y}}

\def\Xbs{\bs{\mc{X}}}

\title{
Nash equilibrium seeking under partial decision information: Monotonicity, smoothness and proximal-point algorithms 
}

\author{Mattia Bianchi, Sergio Grammatico
\thanks{
M. Bianchi and Sergio Grammatico are with the Delft Center for Systems and Control, TU Delft, The Netherlands.
	E-mail addresses: \texttt{\{m.bianchi, s.grammatico\}@tudelft.nl}. This work was partially supported by NWO under research project OMEGA (613.001.702) and by the ERC under research project COSMOS (802348).
}
}

\begin{document}

\maketitle
\thispagestyle{empty}
\pagestyle{empty}

\begin{abstract} 
We address Nash equilibrium problems in a partial-decision information scenario, where each agent can
only exchange information with some neighbors, while its cost function possibly depends on the strategies of all agents. We characterize the relation between several monotonicity and smoothness conditions postulated in the literature. Furthermore, we prove convergence of a preconditioned proximal point algorithm, under a restricted monotonicity property that allows for a non-Lipschitz, non-continuous game mapping. 

\end{abstract}


\section{Introduction}\label{sec:introduction}

\gls{NE} seeking under  partial decision information  has recently attracted considerable research interest, due to its prospect engineering applications as well as theoretical challenges. This scenario arises when, in the absence of a central coordinator, the agents in a network can only rely on the information received from some neighbors, for instance in ad-hoc-networks and sensor positioning problems \cite{DurrStankovicJohanssonACC2011}, \cite{Belgioioso:Nedic:Grammatico;TAC:2021}.  
The technical goal is the distributed computation of a \gls{NE}; the main complication is that the cost function of each agent may depend on the decision variables of  other non-neighboring agents. To cope with the lack of knowledge, each agent estimates and
tries to reconstruct the strategies of all the competitors \cite{TatarenkoShiNedic_Geometric_TAC2021,GadjovPavel_Passivity_TAC2019} 
(or an aggregation value \cite{KoshalNedicShanbhag_Aggregative_OR2016,GadjovPavel_Aggregative_TAC2021})  via peer-to-peer communication. 

In fact, most existing  methods resort to pseudogradient and consensus-type dynamics \cite{YeHu_Consensus_TAC2017,DePersis:Grammatico:Averaging:AUT2019}. Some works studied linearly convergent algorithms, for games without coupling constraints \cite{TatarenkoShiNedic_Geometric_TAC2021,Bianchi_Timevarying_LCSS2021}. Other authors focused on generalized games, for example resorting to an operator-theoretic approach and forward-backward dual methods \cite{GadjovPavel_Aggregative_TAC2021,Pavel_GNE_TAC2020}. All these schemes mainly suffer three drawbacks. 

The first is that gradient-based methods typically require restrictive \emph{monotonicity} assumptions for convergence.  For instance, all the cited works postulate strong monotonicity of the game mapping. Weaker conditions are sometimes sufficient if allowing for vanishing stepsizes: strict monotonicity in the seminal work \cite{KoshalNedicShanbhag_Aggregative_OR2016}, cocoercivity for the 
generalized games in \cite{Belgioioso:Nedic:Grammatico;TAC:2021}. Remarkably, mere monotonicity was recently assumed in \cite{Lei:Shanbhag:Monotone:CDC:2020}, via an additional diminishing Tikhonov regularization. Nonetheless, vanishing stepsizes are undesirable as they affect the convergence speed. Most recently, the authors of \cite{Gadjov:Pavel:Hypomonotone:arXiv:2021} proposed a continuos-time gradient-based method  for (hypo)-monotone games under a novel inverse Lipschitz assumption. 

The second is that the agents' costs must be \emph{differentiable} with \emph{Lipschitz} gradient \cite{YeHu_Consensus_TAC2017,Pavel_GNE_TAC2020}; in turn this ensures that the pseudogradient mapping of the game is Lipschitz. As the game mapping is a global operator, implementing, in a distributed setup, the common alternatives employed in nonsmooth optimization (linesearch or adaptive steps) seems far from trivial. The third is that, due to  partial decision information, the stepsizes must be chosen very small, increasing the number of  iterations for convergence. Importantly, this also translates in prohibitive communication cost, as the agents need to exchange data at each step. 

A possible solution to remedy all three limitations is the proximal-point method \cite[Th.~23.41]{BauschkeCombettes_2017}.
Although a direct implementation in games results in double layer schemes (where the agents have to communicate virtually infinite time between iterations \cite{Scutari:Complex:TIT:2014,Yi:Pavel:Double:TCNS:2019}), in our recent work  \cite{Bianchi:GNEPPP:AUT:2022,Bianchi:PPP:CDC:2020} we have shown that an efficient method can be obtained via \emph{preconditioning} --for the case of games with strongly monotone and Lipschitz mapping. The result is that, at the price of some additional \emph{local} complexity, the number of iterations and communications for convergence to a \gls{NE} can be substantially reduced.  

In this paper we further leverage the properties of \glspl{PPA} to deal with the other two  issues: monotonicity and smoothness. Our  contributions are summarized as follows:
\begin{itemize}
    \item We compare a significant group of monotonicity and smoothness assumptions employed  in the partial decision information literature. We characterize the relations between the conditions, and exemplify their restrictiveness (§\ref{sec:taxonomy});
    \item We prove convergence of our fully distributed \gls{NE} seeking \gls{PPP} algorithm, under the restricted monotonicity of an augmented operator. Our condition is remarkably weaker than that recently proposed in \cite[Th.~2]{Huang:Hu:DR:CDC:2021}  (for a  Douglas-Rachford algorithm). In particular, we
   do not assume  strong monotonicity, nor continuity of the game mapping --which requires a different limiting argument compared to \cite[Th.~2]{Bianchi:GNEPPP:AUT:2022}. Interestingly,  nonsmoothness only affects the \emph{local} optimization problems of the agents (§\ref{sec:PPP}).
\end{itemize}
The proofs are in Appendix. 
\section{Preliminaries}
\subsubsection{Notation} $[A]_{i,j}$ is the element on  row $i$ and column $j$ of a matrix $A$.  $\otimes$ denotes the Kronecker product. $I_n$ is an identity matrix of dimension $n$; we may omit the subscript if there is no ambiguity. 

\subsubsection{Euclidean spaces} Given a positive definite matrix $ \R^{q\times q} \ni P \succ 0$,  $\mc{H}_P\coloneqq (\R^q,\langle \cdot , \cdot \rangle _P)$ is the Euclidean space obtained by endowing $\R^q$ with the $P$-weighted inner product  $\langle x , y \rangle _P=x^\top P y$, and $\| \cdot \|_P$ is the associated norm; we omit the subscripts if $P=I$. Unless otherwise stated, we always assume to work in $\H=\H_I$. 

\subsubsection{Operator-theoretic background \cite{BauschkeCombettes_2017}}
	A set-valued  operator $\mc{F}:\R^q\rightrightarrows \R^q$ is characterized by its graph
	$\gra (\mc{F})\coloneqq \{(x,u) \mid u\in \mc{F}(x)\}$. $\dom(\mc{F})\coloneqq \{x\in\R^q| \mc{F}(x)\neq \varnothing \}$,
	$\fix\left( \mathcal{F}\right) \coloneqq  \left\{ x \in \R^q \mid x \in \mathcal{F}(x) \right\}$ and $\zer\left( \mathcal{F}\right) \coloneqq  \left\{ x \in \R^q \mid 0 \in \mathcal{F}(x) \right\}$ are the domain, set of fixed points and set of zeros, respectively. $\mc{F}^{-1} $ denotes the inverse operator of $\mc{F}$, defined as $\gra (\mc{F}^{-1})=\{(u,x)\mid (x,u)\in \gra(\mc{F})\}$.
	$\mathcal{F}$ is
	(strictly, $\mu$-strongly, $\nu$-hypo-) monotone in $\mc{H}_P$ if $\langle u-v , x-y\rangle_P \geq 0$ ($>0$, $\geq \mu\|x-y\|^{2}_P$, $\geq - \nu \|x-y\|^{2}_P$) for all $(x,u)$,$(y,v)\in\gra(\mc{F})$; 
	 we omit the indication ``in $\mc{H}_P$'' whenever $P=I$. $\mc{F}$ is maximally monotone if it is monotone and there is no operator $\mc{A}$ such that $\gra(\mc{F}) \subset \gra(\mc{A})$. A single-valued operator $\Fmc:\R^q\rightarrow \R^q$ is ${\beta}$-cocoercive if $\langle x-y,\Fmc(x)-\Fmc(y) \geq \beta \| \Fmc(x)-\Fmc(y)\|^2 $ for all $x,y\in\R^q$ (equivalently, $\Fmc^{-1}$ is $\beta$-strongly monotone); is $R$-inverse Lipschitz if $R \|\Fmc{x}-\Fmc{y}\rangle \| \geq \|x-y\|$ (equivalently, $\Fmc^{-1}$ is $R$-Lipschitz). 
	$\Id$ is the identity operator. ${\rm J}_{\mathcal{F} }\coloneqq (\Id + \mathcal{F} )^{-1}$ denotes the resolvent operator of $\mathcal{F} $. For a function $\psi: \R^q \rightarrow \R \cup \{\infty\}$, $\dom(\psi) \coloneqq  \{x \in \R^q \mid \psi(x) < \infty\}$;  its subdifferential operator is
	$\partial \psi: \dom(\psi) \rightrightarrows \R^q:x\mapsto  \{ v \in \R^q \mid \psi(z) \geq \psi(x) + \langle v \mid z-x \rangle , \forall  z \in {\rm dom}(\psi) \}$; if $\psi$ is differentiable and convex, $\partial \psi=\nabla \psi$. For a set  $S\subseteq \R^q$, $\iota_{S}:\R^q \rightarrow \{ 0, \infty \}$ is the indicator function, i.e., $\iota_{S}(x) = 0$ if $x \in S$, $\infty$ otherwise; $\nc_{S}: S \rightrightarrows \R^q:x\mapsto \{ v \in \R^q \mid \sup_{z \in S} \, \langle v \mid z-x \rangle \leq 0  \}$ is the normal cone operator of $S$. If $S$ is closed and convex, then $\partial \iota_S=\nc_S$
	 and $(\Id+\nc_S)^{-1}=\proj_S$ is the Euclidean projection onto  $S$.
Given $\mc{F}:S\rightarrow \R^q$, the variational inequality VI$(\mc{F},S)$ is the problem of finding $x^\star\in S$ such that $\langle \mc{F}(x^\star)\mid x-x^\star\rangle \geq 0$, for all $x \in S$ (or, equivalently, $x^\star$ such that $\0\in\mc{F}(x^\star)+\nc_S(x^\star)$).

\begin{definition}[Restricted monotonicity] \label{def:restricted} An operator $\Fmc:\R^q \rightrightarrows \R^q$ is restricted (strictly, $\mu$-strongly) monotone in $\H_P$ with respect to a set $\Sigma \neq \varnothing$ if $\langle x - x^\star , u- u^\star \rangle _P   \geq 0 $ ($>0$, $\geq \mu \|\omega-\omega^\star\|^2_P $) for all $(x,u) \in \gra(\Fmc)$, $(x^\star,u^\star) \in \gra(\Fmc)$ with $x^\star\in\Sigma$. We omit the characterization in ``$\H_P$'' whenever $P=I$.
\end{definition}

This definition slightly generalizes that in \cite[Def.~1]{Bianchi:GNEPPP:AUT:2022}, which only consider the zero set; note that $\Fmc$ is allowed to be set-valued on  $\omega^\star \in \Sigma$.

\emph{Proximal point algorithm}: For an operator $\Fmc:\R^q\rightrightarrows \R^q$ with $\zer(\Fmc) \neq \varnothing$, we consider the problem of finding a point $x^\star\in \zer(\Fmc)$. The iteration 
\begin{align}\label{eq:ppa}
(\forall k\in\N) \quad  x^{k+1} \in \res_\Fmc (x^k) = (\Id+\Fmc)^{-1}x^k 
\end{align}
is called \gls{PPA}. Note that at each iteration \eqref{eq:ppa} involves solving for $x^{k+1}$ the regularized inclusion $\0 \in \Fmc(x^{k+1}) +  x^{k+1}- x^{k}$. By definition, $\fix(\res_\Fmc) =\zer(\Fmc)$. If $\Fmc$ is maximally monotone, then  $\res_\Fmc$ is single valued and $\dom(\res_\Fmc)=\R^q$, so \eqref{eq:ppa} is uniquely defined; moreover, $x^k$ converges to a point in $\zer(\Fmc)$.



\section{Mathematical setup}\label{sec:prob}

\subsection{The game}
 Let $ \mc I:=\{ 1,\ldots,N \}$ be a set of agents, where each agent $i\in \mc{I}$ chooses its strategy (i.e., decision variable) $x_i$ 
 from its local decision set $\textstyle \X_i \subseteq \R^{
n_i}$.  We denote by $x := \col( (x_i)_{i \in \mc I})  \in \X $  the stacked vector of all the agents' strategies, with $\ \X := \X_1\times\dots\times\X\subseteq \R^n$ the overall decision space and $n:=\textstyle \sum_{i\in \mc{I}} n_i$.
Agent $i \in \mc I$ aims to minimize an objective function $f_i(x_i,x_{-i})$, depending both on the local variable $x_i$ and on the strategies of the other agents $x_{-i}:= \col( (x_j)_{j\in \mc I\backslash \{ i \} })$.
%
The game is  represented by $N$ inter-dependent optimization problems
\begin{align} \label{eq:game}
\forall i \in \mc{I}:
\quad \underset{y_i \in \Omega_i}{\argmin}   \; f_i(y_i,x_{-i}).
\end{align} 
The mathematical problem we consider  is the distributed computation of a \gls{NE},  a set of strategies simultaneously solving all the problems in \eqref{eq:game}.

\begin{definition}
	A Nash equilibrium is a set of strategies $x^{*}=\operatorname{col}\left((x_{i}^{*}\right)_{i \in \mathcal{I}})$ such that, for all $i \in \mc{I}$,
	\[
	\forall i \in \mathcal{I}: \quad x_{i}^{*} \in \underset{y_i\in \Omega_i}{\argmin} \, f_{i}\left(y_{i}, x_{-i}^{*}\right).
	\]
\end{definition}
	\smallskip
Throughout, we restrict our attention to convex games. The following are standard regularity conditions.  
\begin{assumption}[Convexity]\label{asm:convexity}
For each $\i$, the set $\mc{X}_i$ is nonempty, closed and convex; the fuction $f_i$ is continuous and the function $f_i(\cdot,x_{-i})$ is convex for any $x_{-i}$. 
\end{assumption}

Furthermore, we assume existence of a solution. 
\begin{assumption}[Existence]\label{asm:existence}
The game \eqref{eq:game} admits at least one Nash equilibrium. 
\end{assumption}

Sufficient conditions for existence of a \gls{NE} (e.g., compactness of $\X$) can be found, for instance, in \cite{FacchineiPang_PalomarEldar_2009}. 

\subsection{The communication Network}
The agents can exchange information with some neighbors
over an undirected  communication network $\mathcal G(\mc{I},\mc{E})$. The unordered pair $(i,j) $ belongs to the set of edges $\mc{E}$ if and only if agent $i$ and $j$ can mutually exchange information. 
We denote: $W\in \R^{N\times N}$ the  weight matrix of $\mc{G}$, with $w_{i,j}:=[W]_{i,j}$ and  $w_{i,j}>0$ if $(i,j)\in \mc{E}$, $w_{i,j}=0$ otherwise;
$\mc{N}_i=\{j\mid (i,j)\in \mc{E}\}$ the set of neighbors of agent $i$. 

\begin{assumption}[Connectivity]
\label{asm:connectedness}
	The communication graph $\mathcal G (\mc{I},\mc{E}) $ is undirected and connected. 
	The weight matrix $W$ satisfies the following conditions:
	\begin{itemize}[topsep=0em]
	\item[(i)] \emph{Symmetry: }$W=W^\top$;
		\item[(ii)] \emph{Self loops: }$w_{i,i}>0$ for all $i\in \mc{I}$;
		\item[(iii)] \emph{Double stochasticity: }$W\1_N=\1_N, \1^\top W=\1^\top$.
\end{itemize}
We denote by $\sigma \coloneqq \sigma_{N-1}(W)<1 $ the second largest singular value of $W$.
\end{assumption}

The requirements (ii)-(iii) in \cref{asm:connectedness} are intended to ease the notation and they are not strictly necessary; these conditions can for example  be satisfied by assigning Metropolis weights \cite[§2]{Bianchi_Timevarying_LCSS2021}. 

\subsection{The partial decision information scenario}

We consider the so-called partial decision information setup, where  agent $i\in\mc{I}$ can only access its own feasible set $\X_i$  and an analytic expression of its private cost $f_i$,  but cannot access the strategies of all the competitors $x_{-i}$.  Therefore, each agent $i$ is unable to evaluate the actual value of $f_i(x_i,x_{-i})$. 
Instead, each agent keeps an estimate of all other agents' actions \cite{GadjovPavel_Passivity_TAC2019}, \cite{KoshalNedicShanbhag_Aggregative_OR2016}, \cite{DePersis:Grammatico:Averaging:AUT2019}, and aims at reconstructing the actual values, only based   information exchanged locally with neighbors over the communication graph $\mc{G}$. We denote $\x_{i}=\operatorname{col}((\x_{i,j})_{j\in \mc{I}})\in \R^{n}$,  where $\x_{i,i}:=x_i$ and $\x_{i,j}$ is agent $i$'s estimate of agent $j$'s strategy, for all $j\neq i$; $\x_{j,-i}=\col((\x_{j,l})_{l\in\mc{I}\backslash \{ i \}})$; $\x=\col((\x_i)_{i\in\mc{I}})\in\R^{Nn}$ the overall estimate vector; $\x_{-i} =\col((x_j)_{j\in\I\backslash\{i\}}) $. Let
	\begin{align}\label{eq:Rmc}
\mathcal{R}_{i}:=&\left[ \begin{array}{lll}{{0}_{n_{i} \times n_{<i}}} & {I_{n_{i}}} & {\0_{n_{i} \times n_{>i}}}\end{array}\right], 
\end{align}
where $n_{<i}:=\sum_{j<i,j \in \mathcal{I}} n_{j}$, $n_{>i}:=\sum_{j>i, j \in \mathcal{I}} n_{j}$. 
In simple terms,  
$\mathcal R _i$ selects the $i$-th $n_i$-dimensional component from an $n$-dimensional vector, i.e.,  $\mathcal{R}_{i} \x_{i}=\x_{i,i}=x_i$.
Let also $\mathcal{R}:=\operatorname{diag}\left((\mathcal{R}_{i})_{i \in \mathcal{I}}\right)$,
 so that $x=\mathcal{R} \x$.
 
\subsection{Game mapping, extended mapping, augmented operators}

Under \cref{asm:convexity}, a strategy $x^\star$ is a \gls{NE} of the game \eqref{eq:game} if and only if 
\begin{align} \label{eq:NEinclusion}
{\0_{n}} & \in F\left(x^{*}\right)+\mathrm{N}_{\X}\left(x^{*}\right),
\end{align}
where $F:\R^n \rightrightarrows \R^n$ is the \emph{game mapping}
\begin{align}
\label{eq:pseudo-gradient}
F(x):=\operatorname{col}\left( (\subd{x_i} f_i(x_i,x_{-i}))_{i\in\mathcal{I}}\right)
\end{align}
(in fact, \eqref{eq:NEinclusion} are the first order optimality conditions of each convex problems in \eqref{eq:game}). Typically, distributed \gls{NE} seeking methods works under some monotonicity assumption on $F$. Since we deal with the partial decision information scenario, it is also useful to introduce the \emph{extended game mapping} 
\begin{align}\label{eq:extended_pseudogradient}
\Fbs(\x)\coloneqq \operatorname{col}\left( (\subd{x_i} f_i(x_i,\x_{i,-i}))_{i\in\mathcal{I}}\right)
\end{align}
where the subdifferentials are computed on the \emph{estimates}, and the \emph{extended} operators
\begin{align}
\label{eq:Faop}
    \Fa (\x) & \coloneqq \alpha \Fop (\x) + (I_{Nn} - \Wbs) \x
    \\
    \A (\x) & \coloneqq \Fa (\x) \label{eq:Aop} +\nc_{\Xbs}(\x),
\end{align}
where $\alpha>0$ is a design parameter, $\Wbs \coloneqq W\otimes I_n$, $\Xbs \coloneqq \{ \x\in\R^{Nn} \mid \Rmc \x \in \X\}$. The following well-known result (e.g., \cite[Prop.~1]{TatarenkoShiNedic_Geometric_TAC2021}) provides an extension of the inclusion \eqref{eq:NEinclusion} to the estimate space. 

\begin{lemma}\label{lem:VIequivalence}
	The following statements are equivalent:
	\begin{itemize}[topsep=0em] 
		\item [i)] $\bs{x^\star}=\1_N\otimes x^\star$, with  $x^\star\in\X$  a \gls{NE} of the game  \eqref{eq:game};
		\item[ii)] $\0_{Nn}\in \A(\x^\star)$. 
	\end{itemize}
\end{lemma}

In particular, \cref{asm:existence}  implies that $\zer(\A) \neq \varnothing$. 


\section{Towards a  taxonomy
of assumptions}\label{sec:taxonomy}

In recent years, distributed \gls{NE} seeking under partial decision information has been studied under a variety of conditions on the operators $F, \Fop, \Fa, \A$.  Some of the assumptions postulated have not been exemplified, nor it is evident how restrictive they are --in theory and in practice. Towards a solution of this issue, we start by considering the following, representative, conditions. 

\begin{condition}\label{C1}
The operator $\Fop$ is maximally monotone. 
\end{condition}

\begin{condition}\label{C2}
The operator $\Fop$ is restricted monotone with respect to $\zer(\A)$. 
\end{condition}

\begin{condition}\label{C3}
There exists $\alpha \geq 0$ such that the operator $\Fa$ is maximally monotone.
\end{condition}

\begin{condition}\label{C4}
There exists $\alpha \geq 0$ such that the operator $\Fa$ is restricted monotone with respect to $\zer(\A)$.
\end{condition}

\begin{condition}\label{C5}
The  operator $F$ is $\mu$-restricted strongly monotone with respect to $\zer(\A)$ and $\ell$-Lipschitz, for some $\mu > 0$, $\ell>0$. 
\end{condition}

\begin{condition}\label{C6}
The  operator $F$ is $\mu$-strongly monotone and $\ell$-Lipschitz, for some $ \mu > 0$, $\ell>0$. 
\end{condition}

\begin{condition}\label{C7}
The operator $F$ is $\nu$-hypomonotone, $\ell$-Lipschitz, and $R$-inverse Lipschitz, for some $\nu\geq 0$, $\ell>0$, $R>0$, $R\nu <1$. 
\end{condition}

\begin{condition}\label{C8}
The operator $F$ is strictly monotone and $\ell$-Lipschitz, for some $\ell > 0$. 
\end{condition}

\begin{condition}\label{C9}
The operator $F$ is $\frac{1}{\ell}$ cocoercive for some $\ell>0$. 
\end{condition}

\begin{condition}\label{C10}
The operator $F$ is monotone and $\ell$-Lipschitz, for some $\ell > 0$. 
\end{condition}

Although \condref{C6} is  the most common technical assumption, all these conditions have been formulated in the literature (see \cref{tab:1}), except for \condref{C2} (which is a natural relaxations of \condref{C1}) and \condref{C4} (which we will use to show convergence of our algorithm). The following result characterizes the relation between them. 

\begin{figure}
    \centering
    \includegraphics[width=\columnwidth]{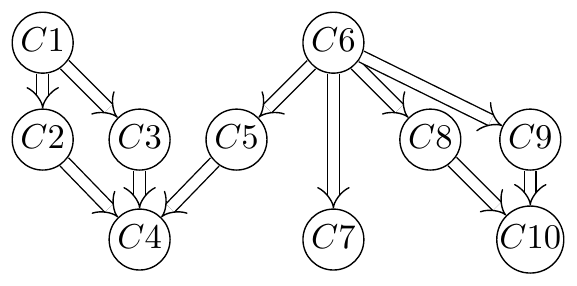}
    \caption{Relations between technical assumptions in monotone games under partial decision information.}
    \label{fig:graph}
\end{figure}

\begin{proposition}\label{prop:relations}
The implications in \cref{fig:graph} hold true. 
\end{proposition}

It can be also shown by counter examples that no other implication exists between the conditions in \condref{C1}-\condref{C10}.

\begin{table}[t]
	\centering
	\newlength{\mywidth}
	\setlength{\mywidth}{0.8em}
	\begin{tabular}{  c c c c  }
		\toprule
		&
		Refs
		& 
		Extra asm.
		& 
		Stepsizes
		\\
		\midrule
		C1 & \cite{GadjovPavel_Passivity_TAC2019,Gadjov:Pavel:Hypomonotone:arXiv:2021} & & Continuous time
		\\
		C3 & \cite{Huang:Hu:DR:CDC:2021,Huang:Hu:Stochastic:arXiv:2022} & & Fixed  \\
		C5 & \cite{TatarenkoShiNedic_Geometric_TAC2021} & &  Fixed \\
		C6 & \cite{Bianchi:GNEPPP:AUT:2022,Pavel_GNE_TAC2020,DePersis:Grammatico:Averaging:AUT2019} & & Fixed \\
		C7 &\cite{Gadjov:Pavel:Hypomonotone:arXiv:2021} &$\X = \R^n$  & Continuous time\\
		C8 & \cite{KoshalNedicShanbhag_Aggregative_OR2016} &  $\X$ compact & Vanishing\\
		C9  & \cite{Belgioioso:Nedic:Grammatico;TAC:2021} & $\X$ compact & Vanishing \\
		C10 & \cite{Lei:Shanbhag:Monotone:CDC:2020} &  $\X$ compact  & Vanishing \\
		\bottomrule
	\end{tabular}
	\vspace{0.3em}
	\caption{\label{tab:1} Technical assumptions in the literature.}
\end{table}

\subsection{ Conditions on the extended pseudogradient}\label{sec:Fbs}

We next prove, under the commonly used assumption that $F$ is single-valued, that \condref{C1} is very restrictive. 

\begin{proposition}[\condref{C1} is trivial]\label{prop:C1istrivial} Assume that $F$ is single valued and continuous. Then, condition \condref{C1} holds if and only if $\grad{} f_i (\cdot, x_{-i})$ is independent of $x_{-i}$, for all $\i$.
\end{proposition}

As the actions $x_{-i}$ are not affecting the optimization problem of agent $i$ (beside possibly for a separable component), there appear to be no reason for agent $i$ to keep estimates (hence, for a partial decision information setup). 

\begin{example}\label{ex:1}
The game defined by $N=2$, $n=2$, $\X = \R^n$, $f_1(x) = {x_1}^2 ({x_2}^2+1)$, $f_2(x) = {x_2}^2({x_1}^2+1)$ has a unique \gls{NE} in $\0$ and satisfies \condref{C2}, but not \condref{C1}.
\end{example}

Although  $\grad{x_i}f_i$ depends on $x_{-i}$ in  \cref{ex:1}, 
the next lemma shows that \condref{C2} is also not of particular interest. 

\begin{proposition}[\condref{C2} is trivial]\label{prop:C2istrivial} Assume that $F$ is single valued and continuous. Then, condition \condref{C2} holds if and only if $\grad{x_i} f_i (x_i^\star,x_{-i})$ is independent of $x_{-i}$, for all $\i$, for any $x^\star = (x_i^\star,x_{-i}^\star)$ \gls{NE} of the game \eqref{eq:game}.
\end{proposition}

In particular, \cref{prop:C2istrivial} implies that 
$
    0  \leq\grad{x_i} f_i(x_i^\star,x_{-i}^\star), x_i - x_i^{\star} \rangle   = \langle  \grad{x_i}  f_i(x_i^\star,x_{-i}), x_i - x_i^{\star} \rangle \
$
where the inequality is the first order optimality condition (as $x_i^\star$ solves \eqref{eq:game}). This  means that, for $x_i^\star$ is optimal for agent $i$ regardless of $x_{-i};$ in other terms, \condref{C2} implies that the Nash equilibria are uniquely composed by dominant strategies (as in \cref{ex:1}). This is also a trivial case, as the agents do not need to communicate to compute a \gls{NE}. 
Although the condition in \cref{prop:C2istrivial} might be violated if $F$ is not continuous, this can only happen at discontinuity points, which is quite a pathological case.

\subsection{Conditions on the game primitives}\label{sec:maxmonFa}

\condref{C3}, \condref{C5}-\condref{C10} are directly postulated on the game mapping $F$ and are  the most well-investigated (e.g., they are easy to check if $F$ is a linear operator \cite{BauschkeCombettes_2017,FacchineiPang_VIs_Springer2004,Gadjov:Pavel:Hypomonotone:arXiv:2021}). \condref{C3}, \condref{C5}. \condref{C8} and the recently proposed \condref{C7} imply uniqueness of the equilibrium; methods with linear convergence were proposed under \condref{C6} \cite{TatarenkoShiNedic_Geometric_TAC2021,Bianchi_Timevarying_LCSS2021}, but not \condref{C7}, \condref{C8}. Although \condref{C5} is weaker than \condref{C6} in theory, it is difficult to check without knowledge of the solutions; we have included it because it causes very limited complications in convergence analysis with respect to \condref{C6}: both conditions actually imply  
 that $\Fa$ is Lipschtz and  restricted strongly-monotone with respect to the whole consensus subspace $\bs{E}:=\{\bs{y}\in\R^{Nn} \mid \y=\1_N\otimes y, \ y\in\R^n \}\supset \zer(\A)$ \cite[Lem.~3]{Pavel_GNE_TAC2020}, a much stronger condition that \condref{C4}).
 \condref{C10} and \condref{C9} allow for multiple \glspl{NE}; yet --as for \condref{C8}-- the related methods require not only compact feasible sets (possibly reasonable in practice) but also vanishing steps, which affects the convergence speed.


\subsection{Conditions on the augmented operator}\label{sec:meremonotonicity}

\condref{C3} and \condref{C4} are more abstract and often replaced by more easily checked  sufficient conditions.
For example, 
restricted monotonicity of $\Fa$ with respect to $\bs{E}$   can be checked without knowledge of the solutions, and implies \condref{C4}. 

Despite this complication, \condref{C3} and \condref{C4}
are of interest for nonsmooth games, as we exemplify next. 

\begin{example}\label{ex:2}
Consider the game defined by $N=2$, $n=2$, $\X = \R^n$, $F(x) = \bar F(x)+ \hat F(x)$, with $\bar F(x) = \col({x_1}^3,0)$ and $\hat F(x) = \left[ \begin{smallmatrix} 2 & 1 \\ 1 & 2
\end{smallmatrix}\right] x+\left[ 
\begin{smallmatrix}
5 \\ 4
\end{smallmatrix}\right]
$. As $\bar F$ is monotone and $\hat F$ is strongly monotone, the game admits a unique equilibrium. Conditions \condref{C5}-\condref{C10} are violated, as they  require Lipschitz continuity of $F$; \condref{C2} also fails (as the best response of agent $2$ is $-0.5x_1$ and by \cref{prop:C2istrivial}). However, \condref{C4} holds: to show this, consider the components of the extended game mapping $\bar \Fbs$ and $\hat \Fbs$ corresponding to $\bar F$ and $\hat F$; $\Rmc^\top \bar \Fbs$ is monotone, while 
$ \alpha \Rmc^\top \hat \Fbs+ (I-\Wbs)$ can be made restricted monotone with respect to the consensus subspace by choosing $\alpha >0$ small enough. We can check numerically that \condref{C3} also holds for some $W$, although there is no analytical test available.
\end{example}

\begin{example}\label{ex:3}
Consider \cref{ex:2} but with $\bar F(x) = \col({x_1}^3({x_2}^4+1),0)$ and $\hat F(x) = \left[ \begin{smallmatrix} 2 & 1 \\ 1 & 2  \end{smallmatrix} \right]x$. The game admits a \gls{NE} $x^\star =0$. As $F$ is restricted strongly monotone with respect to $x^\star$, the equilibrium must be unique. As for \cref{ex:2}, it is easy to check that \condref{C4} holds; yet it can be proven that \condref{C3} does not.  
\end{example}

\begin{example}\label{ex:4}
Consider the game defined by $N=2$, $n=2$, $\X = \R^n$, $f_1(x) = {x_1}^2 -|x_1||x_2|$, $f_2(x) = {x_2}^2 +x_2x_1$, where $|\cdot|$ denotes the absolute value. The game admits a unique \gls{NE} in $\0$; moreover, the operator $F$ is set valued, as $f_1$ is not differentiable in the local variable. Nonetheless, it can be  checked that \condref{C4} holds. 
\end{example}


\section{The \gls{PPP} algorithm}\label{sec:PPP}

In this section we consider the fully-distributed proximal-point \gls{NE} seeking method shown in \cref{algo:1}. The iteration  coincides with that studied in \cite{Bianchi:PPP:CDC:2020}, although the terms have been rearranged. The algorithm includes a consensus phase, where the agents exchange and mix their variable vectors. The local actions are then updated according to a proximal-best response with stepsize $\alpha^k>0$ --importantly, the cost function of each agent $i$ evaluated in the estimates $\x_{i,-i}$, and not on the real competitor's actions $x_{-i}$. Note that the algorithm is always well (uniquely) defined, as the update of $x_i$ is the $\argmin$ of a strongly convex function (by convexity of $f_i(\cdot,x_{-i})$ in \cref{asm:convexity}). 

\begin{algorithm}[t] \caption{Fully-distributed \gls{PPP} algorithm } \label{algo:1}
	\vspace{-0.4em}
	\begin{align*}
		\tilde{\x}_{i}^{k} & = \textstyle \frac{1}{2}( \x_{i}^k+\sum_{j=1 }^{N}w_{i,j}\x_{j}^{k})
		\\
		\x_{i,-i}^{k+1} & = \tilde{\x}_{i,-i}^{k}
		\\
		{x}_{i}^{k+1}     & =
		\underset{y \in \Omega_i} {\argmin} \ \left(
		f_i(y,\tilde{\x}_{i,-i}^{k})+\textstyle \frac{1}{\alpha}  \|  y- \tilde \x_{i,i}^k \|^2 \right)
	\end{align*}
\end{algorithm}

\cref{algo:1} can be formulated as a proximal point method applied to the operator $\A$. However, the computation of $(\Id+\A)^{-1}$ cannot be performed in a distributed way (more precisely, it would require the collaborative solution of a regularized game at each iteration, resulting in a scheme with nested layers of communication, see \cite{Scutari:Complex:TIT:2014}). We have shown in \cite{Bianchi:PPP:CDC:2020,Bianchi:GNEPPP:AUT:2022} that this complication can be tackled by preconditioning the operator $\A$ with a preconditioning matrix 
\begin{align}\label{eq:Phi}
    \Phi = I_{Nn}+ \Wbs.
\end{align}
\begin{lemma}[{\cite[Lem.~2]{Bianchi:PPP:CDC:2020}}]\label{lem:ref_prec}
\cref{algo:1} can be written as
\begin{align}
    \x^{k+1} = (\Id+\Phi^{-1}\A)^{-1} (\x^k).
\end{align}
\end{lemma}
\smallskip
This operator-theoretic interpretation is very powerful, as it seamlessly allows to study convergence of analogous proximal-best response schemes even in the presence of inexact updates (i.e., the $\argmin$ is only approximated at each iteration), coupling constraint, acceleration terms \cite{Bianchi:GNEPPP:AUT:2022}. It also immediately shows that the fixed points of \cref{algo:1} coincide with $\zer(\A) =\zer(\Phi^{-1}\A)$ (i.e., they are estimates at consensus at a Nash equilibrium). 


The following theorem is the main result of the paper. It extends the convergence results in  \cite[Th.~3]{Bianchi:GNEPPP:AUT:2022}, formulated under \condref{C6}, to the case of restricted monotone --possibly nonsmoooth-- games (\condref{C4}). 


\begin{theorem}\label{th:main}
Let \cref{asm:convexity,asm:existence,asm:connectedness} hold. 
Assume that \condref{C4} holds for some $\alpha>0$.  Then, the sequence $(\x^k)$ generated by \cref{algo:1} converges to a point $\x^\star = \1_N \otimes x^\star$, where $x^\star$ is a Nash equilibrium of the game in \eqref{eq:game}.
\end{theorem}

\begin{remark} In \cite{Bianchi:PPP:CDC:2020} we have proven (linear) convergence of \cref{algo:1} assuming \condref{C6}; under the weaker \condref{C4}, \cref{th:main} leverages the general results for the proximal-point algorithm of restricted (merely) monotone games \cite{Bianchi:GNEPPP:AUT:2022}. With respect to \cite{Bianchi:GNEPPP:AUT:2022} and to the Douglas-Rachford algorithm in \cite{Huang:Hu:DR:CDC:2021}, we use a different limiting argument in our proof, which does not require $F$ to be Lipschitz continuous (or even continuous). The core idea is to show that the operator $\res_{\Phi^{-1}\A}$ is continuous, even if $\A$ is not (nor is maximally monotone). 
For instance, \cref{th:main} can be applied to the games in \cref{ex:2,ex:3,ex:4}, while \cite[Th.~2]{Bianchi:GNEPPP:AUT:2022}, \cite[Th.~3]{Huang:Hu:DR:CDC:2021} cannot. Our examples also show a significant gap between \condref{C4} and the condition \condref{C3}, employed e.g., in \cite[Th.~3]{Huang:Hu:DR:CDC:2021}. \end{remark}

We conclude this section by sketching some technical extensions of our results. 
To start, our arguments in \cref{th:main} can be readily adapted to the algorithms --for generalized games-- studied in \cite{Bianchi:GNEPPP:AUT:2022}, to show convergence even under \condref{C4}. Moreover, our convergence results would hold assuming the definition of restricted monotonicity proposed in    \cite[Def.~1]{Bianchi:GNEPPP:AUT:2022}, slightly less restrictive than our \cref{def:restricted}.
We also note that we assumed monotonicity properties of $F$ (and similarly for the other game operators) to hold over all $\R^n$; however, the conditions can be relaxed to hold only over the feasible set, if the estimates $\x$'s are initialized in $\X^N$ (since the update in \cref{algo:1} guarantees invariance for this set). The costs in \eqref{eq:game} can be modified to include a  more general (discontinuous) proper, convex, closed function $g_i(x_i)$ (besides the indicator function $\iota_{X_i}$), without particular technical complications.  Much more intriguing is the case of discontinuity in the part of the cost coupled with the other agents (i.e., violating \cref{asm:convexity}): although our convergence arguments do not hold in this case, it would be interesting to verify whether \condref{C3} could be satisfied to apply standard \gls{PPA} results. 



\section{Conclusion and outline}\label{sec:Conclusion}
Besides their efficiency, proximal point algorithms have the advantage of only requiring mild monotonicity and smoothness conditions.
We have compared and analyzed, several assumptions in \gls{NE} seeking under partial decision information, and proved the convergence of a fully distributed \gls{PPP} method under one of the weakest. 

Future work should investigate linear rates  in absence of (restricted) strong monotonicity. One promising option is to leverage inverse Lipschitz properties, which can ensure  contractivity of certain resolvents. Proving convergence in merely monotone regime, under fixed step sizes, is also a challenging open problem.


\begin{appendix}

\subsubsection{Proof of \cref{prop:relations}}
$\condref{C1} \Rightarrow \condref{C2}$, $\condref{C3} \Rightarrow \condref{C4}$, $\condref{C6} \Rightarrow \condref{C5}$, $\condref{C6} \Rightarrow \condref{C8}$, $\condref{C8} \Rightarrow \condref{C10}$: By   definition. 

$\condref{C1} \Rightarrow \condref{C3}$: As $(I-W)$ is a positive semidefinite matrix, the operator $I-\Wbs$ is maximally monotone. Hence, for any $\alpha\geq $, $\Fa = \alpha \Fop+(I-\Wbs) $ is the sum of two maximally monotone operators; moreover, $\dom(I-\Wbs)=\R^{Nn}$, so the conclusion follows by \cite[Cor. 25.5]{BauschkeCombettes_2017}.

$\condref{C2} \Rightarrow \condref{C4}$: $\Fa$ is the sum of a restricted monotone operator and a monotone operator, hence restricted monotone. 

$\condref{C5} \Rightarrow \condref{C6}$: See, for instance, \cite[Lem.~3]{Bianchi:GNEPPP:AUT:2022}.

$\condref{C6} \Rightarrow \condref{C7}$: It follows by definition and  \cite[Prop.~3]{Gadjov:Pavel:Hypomonotone:arXiv:2021}.

$\condref{C6} \Rightarrow \condref{C9}$: See e.g. \cite[Prop.~5]{Gadjov:Pavel:Hypomonotone:arXiv:2021}.

$\condref{C9} \Rightarrow \condref{C10}$: It follows by definition of cocoercivity and the Cauchy--Schwartz inequality.   \hfill $\blacksquare$

\subsubsection{Proof of \cref{prop:C2istrivial}} 
``$\Rightarrow$'': For the sake of contradiction, assume that, for some $\i$, there exist $l \in \{1,2, \dots, n_i\}$, $x_i \in \R^{n_i}$ and a pair of vectors $x_{-i}$ and $x_{-i}'$ such that $[\grad{x_i}f_i (x_i,x_{-i})]_l < [\grad{x_i}f_i (x_i,x_{-i}')]_l  $. By continuity, there exists $\epsilon>0$ such that $[\grad{x_i}f_i (x_i,x_{-i})]_l < [\grad{x_i}f_i (x_i-\epsilon e_l,x_{-i}')]_l $, where $e_l \in \R^n_i$ is the $l$-th vector of the canonical basis. The  monotonicity in \condref{C1}, applied to pair of  estimate vectors $(\x_i,\x_{-i})$, $(\x_i',\x_{-i})$, for any $\x_{-i}$ and $\x_i = (x_i,x_{-i})$, $\x_i' = (x_i+\epsilon e_l, x_{-i}^\prime)$, gives
\begin{align*}
   0 & \leq  \langle  \grad{x_i}f_i (x_i,x_{-i}) - \grad{x_i}f_i (x_i-\epsilon e_l,x_{-i}'), \epsilon e_l \rangle & \\
   & = \epsilon [\grad{x_i}f_i (x_i,x_{-i}) -  \grad{x_i}f_i (x_i-\epsilon e_l,x_{-i}')]_l & <0   
\end{align*}
 which is a contradiction. Because $x_{-i}$, $x_{-i}'$ are arbitrary, we conclude that, for all $\i$, for all $x_i$, and for all  $x_{-i}$, $x_{-i}'$, $\grad{x_i}f_i (x_i,x_{-i}) = \grad{x_i}f_i (x_i,x_{-i}')$.
 
 ``$\Leftarrow$'': By assumption, for any $\i$, $x_i$, $x_i'$, $x_{-i}$, $x_{-i}'$,
 \begin{align*}
    & \hphantom{{}={}}   \langle \grad{x_i} f_i(x_i,x_{-i}) -\grad{x_i} f_i(x_i',x_{-i}'), x_i - x_i'  \rangle  &   \\ 
     & =
     \langle \grad{x_i}f_i(x_i,x_{-i}') -\grad{x_i} f_i(x_i',x_{-i}'), x_i - x_i& \rangle \geq 0,
 \end{align*}
 where the inequality is convexity of $f_i$ in the first argument (\condref{asm:convexity}). Stacking the inequalities for $\i$ retrieves monotonicity of $\Fop$. \hfill $\blacksquare$

\subsubsection{Proof of \cref{prop:C2istrivial}}
``$\Rightarrow$'': For contradiction, assume that there exist $\i$, $l \in \{1,2, \dots, n_i\}$, $x^\star \in \sol$ and $x_{-i}$ such that $[\grad{x_i}f_i (x_i^\star,x_{-i})]_l < [\grad{x_i}f_i (x_i^\star,x_{-i}^\star)]_l  $. By continuity, there exists $\epsilon>0$ such that $[\grad{x_i}f_i (x_i^\star+\epsilon e_l,x_{-i})]_l < [\grad{x_i}f_i (x_i^\star,x_{-i}^\star)]_l$. Restricted  monotonicity in \condref{C2}, applied to pair of  estimate vectors $(\x_i,\x_{-i})$, $(x^\star,\x_{-i})$, for any $\x_{-i}$ and $\x_i = (x_i^\star+\epsilon e_l,x_{-i})$, gives
\begin{align*}
   0 & \leq  \langle  \grad{x_i}f_i (x_i^\star+\epsilon e_l,x_{-i}) - \grad{x_i}f_i (x_i^\star,x_{-i}^\star),  \epsilon e_l \rangle & \\
   & = \epsilon [\grad{x_i}f_i (x_i,x_{-i}) -  \grad{x_i}f_i (x_i-\epsilon e_l,x_{-i}')]_l & < 0   
\end{align*}
 which is a contradiction. Analogously it can be shown that $[\grad{x_i}f_i (x_i^\star,x_{-i})]_l > [\grad{x_i}f_i (x_i^\star,x_{-i}^\star)]_l  $ leads to a contradiction. Hence $\grad{x_i}f_i (x_i^\star,x_{-i}) = \grad{x_i}f_i (x_i^\star,x_{-i}^\star)$. 

 ``$\Rightarrow$'': For any $\i$, $x_i$, $x_{-i}$, $x^\star\in\sol$, by assumption and convexity,
$ \langle  f_i(x_i,x_{-i}) -\grad{x_i} f_i(x_i^\star,x_{-i}\star), x_i - x_i^\star 
 \rangle  =
     \langle \grad{x_i}f_i(x_i,x_{-i}) -\grad{x_i} f_i(x_i^\star, x_{-i}), x_i - x_i^\star \geq 0.
$ \hfill $\blacksquare$

\subsubsection{Proof of \cref{lem:ref_prec}}
	We have
	\allowdisplaybreaks
	\begin{align*}
	\nonumber
	\hspace{-1.5 em}\x^{k+1} \! & \in (\Id+\Phi^{-1}\A)^{-1}\x^{k}
	\\
	\nonumber
	\hspace{-1.5 em}\iff 
	\0_{Nn}  & \in \x^{k+1}+\Phi^{-1} \Ak\x^{k+1}-\x^{k}
	\\
	\nonumber
	\hspace{-1.5 em}\iff 
	\0_{Nn}   & \in \Phi(\x^{k+1}-\x^{k})+\Ak\x^{k+1} 
	\\
\nonumber
	\hspace{-1.5 em}\iff 
	\0_{Nn} & \in \x^{k+1}+\cancel{\bs{W}\x^{k+1}}-\x^{k}-\bs{W}\x^{k} 	+\x^{k+1}
	\\ 
	\label{eq:stepderivation}
	&
    \quad -\cancel{\bs{W}\x^{k+1}}+\alpha^k\mc{R}^\top \bs{F}(\x^{k+1})+ \mathrm{N}_{\bs{\X}}(\x^{k+1}).  \hspace{-1 em}
	\end{align*}
	The lemma follows by writing componentwise the last inclusion, and by recalling that the zeros of the  subdifferential of a strongly convex function coincide with the unique minimum \cite[Th.~16.3]{BauschkeCombettes_2017}. Note that the preconditioning decouples the updates of agent $i$ from the ``future'' (i.e., at $k+1$) value of $\x_{-i}$, enabling distributed implementation. ~\hfill$\blacksquare$


\subsubsection{Proof of \cref{th:main}}
 We start by auxiliary result.
\begin{lemma}\label{lem:continuity}
The operator $\res_{\Phi^{-1}\A}$ is continuous. 
\end{lemma}
\begin{proof}
By \cref{lem:ref_prec} and the explicit form of $\res_{\Phi^{-1}\A}$ in \cref{algo:1}, we just need to show that the function $h:\tilde{\x_{i}} \mapsto {\argmin}_y ( \tilde f_i(y,\tilde{\x}_{i}) +\iota_{\X_i}(y) ) $, with $\tilde f_i(y,\tilde{\x}_{i}) \coloneqq f_i(y,\tilde{\x}_{i,-i})+ \frac{1}{\alpha}\|y- \tilde \x_{i,i}^k \|^2$,   is continuous (since composition of continuous functions is continuous). Consider any converging (bounded) sequence $\tilde{\x}_{i}^k \rightarrow \tilde{\x}_{i}^*$, and define $x_i^k = \tilde h(\tilde{\x}_{i}^k)$, $x^* = h(\tilde{\x}_{i}^*)$. Note that $\tilde f_i$ is strongly convex, so $(x_i^k)_\k$ must also be bounded. 
Moreover, for any diverging subsequence $\bar{K} = (\bar{k}_1,\bar{k}_2,\dots)\subseteq \N$, $\tilde f_i(x_i^{\bar{k}_n},\tilde{\x}_{i}^{\bar{k}_n})\leq f_i(x_i,\tilde{\x}_{i}^{\bar{k}_n})$ for any $x_i\in \X_i$,
and by continuity of  $f_i$ we conclude that $\tilde f_i(x_i',\tilde{\x}_{i}^*)\leq f_i(x_i,\tilde{\x}_{i}^*)$ for all $x_i \in \X_i$ and any $x_i'$ accumulation point of $(x_i^{\bar{k}_n})_{n\in \N}$. Since the minimizer is
unique by strong convexity, we conclude that $x_i' = x^*$, which also means $x_i^k \rightarrow x^*$.
\end{proof}
 We are now in a position to apply the results on proximal-point algorithm for restricted monotone operators in \cite{Bianchi:GNEPPP:AUT:2022}. First, note that the operator $\A$ is restricted monotone with respect to $\zer({\A})$ (because $\Fa$ is so (by assumption) and by monotonicity of the normal cone  \cite[Th.~20.25]{BauschkeCombettes_2017}), i.e., for all $(\x,\ubs), (x^\star,\ubs^\star)\in\gra(\A)$, with $x^\star \in \zer(\A)$
 \begin{align}
     0 & \leq  \langle \ubs-\ubs^\star, \x-\x^\star \rangle \\
        & =   \langle  \Phi^{-1}\ubs- \Phi^{-1}\ubs^\star, \x-\x^\star  \rangle_{\Phi},
 \end{align}
which  shows that $\Phi^{-1}\A$ is restricted monotone with respect to $\zer(\A)$ in $\H_{\Phi}$. Therefore, by \cref{lem:ref_prec} and by applying \cite[Th.~1(i)]{Bianchi:GNEPPP:AUT:2022}, we  infer that the sequence $(\x^k)$ is bounded, hence it admits at least one cluster point, say $\bar \x$. By \cite[Th.~1(ii)]{Bianchi:GNEPPP:AUT:2022}, $\res_{\Phi^{-1}\A}(\x^k) -\x^k \rightarrow 0$; therefore, by continuity in \cref{lem:continuity}, it must be $\bar \x \in \fix(\res_{\Phi^{-1}\A}) = \zer(\A). $ The conclusion follows by \cite[Th.~1(iii)]{Bianchi:GNEPPP:AUT:2022}. 

\end{appendix}

\bibliographystyle{IEEEtran}
\bibliography{library}

\onecolumn
\section{Auxiliary material: Continuity of the $\operatorname{argmin}$ of strongly convex functions}

\vspace{1em}

\begin{lemma}
Let $f:\R^n\times \R^m \rightarrow \R : (x,y)\mapsto f(x,y)$ be a continuous function, and assume that $f(\cdot,y)$ is $\mu$-strongly convex for any $y \in \R^m$. Let $X\subseteq \R^n$ be convex closed. Then the (single valued, full domain) function 
\begin{align}
    y \mapsto g(y) = \operatorname{argmin}_{x\in X} f(x,y)
\end{align}
is continuous. \hfill $\square$
\end{lemma}

\begin{proof}
For any given sequence $y^k \rightarrow y^*$ (converging, hence bounded), we will show that $x^k := g(y^k) \rightarrow  g(y^*)=: x^*$; this is the definition of continuity of $g$. 

1) First, we show that the sequence $x^k$ is bounded. Let $Y$ be a compact set containing $(y^k)_{k\in\N}$. Let $x_0\in X$ and 
\begin{align}
    l_0 & := \max_{y\in Y} \  f(x_0, y)  \\
    l_1 & := \min_{x\in \partial B(x_0, 1), y\in Y} \ f(x,y) 
\end{align}
where $\partial B( x_0,1) = \{x \in \R^n \mid \|x-x_0\| = 1 \}$ is the boundary of the unit ball centered at $x_0$; the $\min$ and $\max$ must be achieved because the domains are compact. Let $d \in \R^n$ be any unitary vector, i.e., $\|d\| =1$; $x_1 := x_0 +d \in \partial B (x_0,1)$; $x_2 = x_0 + Md$, for some scalar such that $M>1$ and 
\begin{align}
    M > 2 \frac{l_0 - l_1}{\mu}+1.
\end{align} Then, 
\begin{align}
    x_1 = \frac{M-1}{M}x_0 + \frac{1}{M} x_2.
\end{align}

By definition of strong convexity, this means that, for all $y \in Y$
\begin{align}
    l_1 & \leq f(x_1,y) 
    \\
    & \leq  \frac{M-1}{M} f(x_0,y)+ \frac{1}{M}f(x_2,y) - \frac{1}{2}\mu \frac{M-1}{M} \frac{1}{M}\|x_0 - x_2\|^2
    \\
    & =   \frac{M-1}{M} f(x_0,y)+ \frac{1}{M}f(x_2,y) - \frac{1}{2}\mu (M-1)
\end{align}

Assume by contradiction that there exists $y \in Y$ such that $f(x_2,y) \leq f(x_0,y)$. Then, since $f(x_0,y) \leq l_0 $ the previous inequality implies  
$
    l_1 - l_0 \leq - \frac{1}{2}\mu(M-1)
$, which contradicts the assumption on $M$. We conclude that, for any $y \in Y$, for all $x$ such that $\| x_0 -x \| >M $, $f(x_0,y) < f(x,y)$. In turn, this means that for all $y\in Y$, $\| g(y) \| < \| x_0\|+M$, i.e., $g$ is uniformly bounded over $Y$. 

2) Consider any accumulation point $x'$ of $(x^k)$ (one exists by boundedness), and let $\bar{K} = (\bar{k}_1,\bar{k}_2,\dots)\subseteq \N$ be a diverging subsequence such that $x^{\bar k_n} \rightarrow x'$. Since $f(x^{\bar k_n},y^{\bar k _n} ) \leq f(x,y^{\bar k_n}) $ for all $x\in X$, then, by continuity of $f$, $f(x',y^*) \leq f( x, y^*)$ for all $x\in X$. Since the minimizer must be unique by strong convexity, we have $x' = x^*$. In particular, this shows that $x^*$ is the unique accumulation point of $x^k$: therefore, $x^k \rightarrow x^*$. 
\end{proof}

\end{document}